\documentclass[review]{elsarticle}

\usepackage{lineno,hyperref}
\usepackage{xcolor}
\modulolinenumbers[5]

\journal{Arxiv.org}

\biboptions{authoryear}

\usepackage{float}
\usepackage[utf8]{inputenc}
\usepackage{color}
\usepackage{tikz}
\usepackage{amsmath}
\usepackage{amsfonts}
\usepackage[nolist]{acronym}
\usepackage[linesnumbered,noend]{algorithm2e}
\usepackage{float}
\usepackage{subcaption}

\allowdisplaybreaks

\begin{document}
\title{Fast delta evaluation for the Vehicle Routing Problem with Multiple Time Windows}
\author[dtu]{Rune Larsen\corref{cor1}}
\ead{rular@dtu.dk}
\author[dtu]{Dario Pacino}
\ead{darpa@dtu.dk}
\cortext[cor1]{Corresponding author}
\address[dtu]{DTU Management Engineering, Technical University of Denmark}

\begin{abstract}
In many applications of vehicle routing, a set of time windows are feasible for each visit, giving rise to the Vehicle Routing Problem with Multiple Time Windows (VRPMTW). We argue that such disjunctions are problematic for many solution methods, and exemplify this using a state of the art Adaptive Large Neighbourhood Search heuristic.
VRPMTW comes in two variants depending on whether the time used en route must be minimised. A more compact and corrected mathematical formulation for both variants of the problem is presented, and new best solutions for all but six benchmark instances of VRPMTW without time minimisation is found. A new solution representation for VRPMTW with time minimisation is presented, its importance is demonstrated and it is used to find new best solutions for all but one benchmark instance of VRPMTW with time minimisation.
\end{abstract}




\maketitle

\begin{acronym}
\acro{VRP}{the Vehicle Routing Problem}
\acro{VRPTW}{the Vehicle Routing Problem with Time Windows}
\acro{VRPMTW}{Vehicle Routing Problem with Multiple Time Windows}
\acro{TSP}{Travelling Salesman Problem}
\acro{ALNS}{Adaptive Large Neighbourhood Search}
\end{acronym}


\section{Introduction}
The literature on \ac{VRP} is very extensive, both concerning applied methods and on variations of the problem~\citep{Braekers2015TheReview}. A specific variation that has recently caught the attention of the research community is the \ac{VRPMTW}. As mentioned in \cite{Belhaiza2014AWindows}, applications of this problem can be found in delivery operations of furniture and electronic retailers, where a choice of delivery periods might be offered. In countries such as Denmark, deliveries to e.g. large grocery distribution centres require the booking of available time slots on an on-line system. Missing a booked window might result in extreme delays for the entire route. A solution to the \ac{VRPMTW}, where the available time slots are modelled as multiple time windows, will, not only, optimise the route plan but also advice on which time slot to book. As similar case-study can also be found in Portuguese companies~\citep{Amorim2014AStudy}.

The \ac{VRPMTW} is not a new problem in the literature, though the list of works is very small. To the best of the authors' knowledge, the problem was introduced by \cite{Favaretto2007AntVisits}. The authors provided a novel set of benchmark instances and proposed an ant colony optimisation approach for the single and multiple visits variant of the problem. The current state-of-the-art approach for the \ac{VRPMTW} is the hybrid variable neighbourhood tabu search presented by \cite{Belhaiza2014AWindows}. They generalised a set of well-known neighbourhood operators for the VRPMTW: 2-exchange, 2-opt, relocate, 3-node swap, and a 3-exchange. The operators are applied within a tabu search framework where shaking phases and restarts are applied to increase diversification. To assist in time minimisation, they present a recursive approach to calculate the optimal allocation of waiting times along the route, and thus implicitly selecting the time windows used.
Subsequently, the same approach is used within a game theory framework for multiple-criterion optimisation by the same authors~\citep{Belhaiza2016AProblem}. The use of multiple time windows has also been studied in other related problems, from which inspiration and knowledge can be drawn. Applied to the \ac{TSP}, \cite{Pesant1999OnProblem} show the flexibility of constraint programming when modelling e.g. multiple time windows. \cite{Baltz2014ExactSelection} show how simple feasibility checks during insertion operators can handle multiple time windows, while \cite{PaulsenHeuristicWindows} use dynamic programming to evaluate the optimal set of time windows to use when time minimising an entire route in the context of genetic algorithms.
Multiple time windows have also been studied in the team orienting problem~\citep{Souffriau2010TheWindows,Lin2015AWindows,Tricoire2010HeuristicsWindows}, where they have mostly been handled by feasibility checks and not as part of the optimisation. Recently multiple time windows have also been included into a mathematical formulation of the long-haul transportation problem~\citep{Neves-Moreira2016ALocations}.

In this paper, we  present  a  more  compact and a corrected mathematical formulation for the problem, we describe how to adapt the \ac{ALNS} framework~\citep{ALNS2006} for the \ac{VRPMTW} with a novel solution representation allowing delta evaluation to get the cost of insertions when time minimising for the \emph{optimal} allocation of time windows. As opposed to most works on \ac{VRPTW} wherefrom we draw inspiration, our algorithm can be applied to both the distance and time minimisation variant of the problem. Computational results show that the adapted \ac{ALNS} can find better solutions than the current state-of-the-art heuristic while retaining the modelling flexibility of the \ac{ALNS} framework.
We also measure the importance of each component of the solution method.


The remainder of the paper is organised as follows: First, we introduce a compact formulation for the problem (Section~\ref{sec:problem}). Section \ref{sec:ALNS} describes the \ac{ALNS} used in this paper, focusing on the fast evaluation of the insertions. Section~\ref{sec:config} documents the setup of the \ac{ALNS} framework to solve the \ac{VRPMTW}, and quantifies the importance of each component. We then demonstrate the effectiveness of the proposed heuristic on the instances from the literature in Section~\ref{sec:results}, and finally, we conclude on our findings in Section~\ref{sec:conclusions}.

\section{Problem definition and formulation}
\label{sec:problem}
We consider two variants of the VRPMTW. In both versions, we associate a cost with the travel distance and the number of vehicles used. In the second version, we are also minimising the time used by the vehicles.

The formulation is a simplified version of the model from \citep{Belhaiza2014AWindows}. We use fewer variables, but model the same problem. Eliminated variables include the flow carried along each arc $y^k_{ij}$, $r^k$ indicating if a vehicle is used, $w^k_i$ representing the waiting time at each customer, $q_d^k$ representing the amount loaded at depot, $d_k$ representing the duration of a vehicles route, and finally $b_i^k$ representing the departure time each vehicle at each customer. Furthermore, $v^p_i$ and $z_i^k$ have been merged into a new $z^k_{ij}$, indicating if a visit is served by a vehicle in a given time window. The number of constraint types are reduced from 20 in \citep{Belhaiza2014AWindows} to 13, easing readability further.
Lastly, we adjusted the objective to reflect the one used in the evaluation of the instances in both papers. Namely the service time $s_{i}$ has been added when minimising time, which is currently missing in \citep{Belhaiza2014AWindows}.

 Assuming that an unlimited set of vehicles $K$ is available, the problem is formulated on a graph $G=(V,A)$, in which the set of vertices are partitioned into $V=N\cup \{d_{start},d_{end}\}$ with $N$ being a set of vertices representing the visits, $d_{start}$ being the start depot, and $d_{end}$ being the end depot for the vehicles. In the mathematical model we assume an heterogeneous fleet to match \citep{Belhaiza2014AWindows}, where each vehicle $k\in K$ has a capacity $Q_k$ and a maximum route duration $T^k$. For the remainder of the paper we operate with a homogeneous fleet, as that is the case for the test instance. Each visit $i\in N$ has a cargo demand $q_i$ and a service time $s_i$. A vehicle is allowed to service a visit only within a fixed set of time-windows. Let $P_i$ be the set of time-windows for visit $i$. Each time-window $p\in P_i$ is then defined within the interval $[l^p_i,u^p_i]$. The travel time between two visits $i,j\in V$ is denoted $t_{ij}$, to which we associate a cost $c_{ij}^k$ if travelled by vehicle $k$. An additional cost $c_k$ is also paid for using vehicle $k$. An overview of the notation can be found in Table~\ref{tab:notation}.
\begin{table}[h]
\begin{center}
\begin{tabular}[h]{ l l }
\hline
Parameter                 & Description \\ \hline
$N$                       & The set of vertices representing the visits\\
$V$                       & The set of vertices including the depots\\
$P_i$                     & The set of time-windows for visit $i\in V$ and the time horizon for the depots\\
$t_{ij}\in \mathbb{R}$    & The driving time between visit $i,j\in V$ \\
$s_i\in \mathbb{R}$       & The service time at visit or depot $i\in V$ \\
$c^k\in \mathbb{R}$       & The cost of using vehicle $k\in K$ \\
$c^k_{ij}\in \mathbb{R}$  & The cost of using vehicle $k\in K$ between visits $i,j \in V$ \\
$q_i\in \mathbb{R}$       & Demand of visit $i\in N$\\
$Q_k\in \mathbb{R}$       & Capacity of vehicle $k\in K$ \\
$l^p_i\in \mathbb{R}$     & Lower bound of time window $p\in P_i$ of visit $i\in N$\\
$u^p_i\in \mathbb{R}$     & Upper bound of time window $p\in P_i$ of visit $i\in N$\\
$T^k\in \mathbb{R}$       & Maximum route duration    for vehicle $k\in K$\\
$B\in \mathbb{B}$         & 1 if and only if time usage is penalised \\
$M\in \mathbb{R}$         & A large constant \\
\hline
\hline Variable & Description \\ \hline 
$x^k_{ij} \in \mathbb{B}$ & 1 if and only if vehicle $k\in K$ travels between visits $i,j \in V$ \\
$w^k_i\in \mathbb{R}$     & Waiting time for vehicle $k\in K$ at visit $i\in N$ \\
$z^k_{ip}\in \mathbb{B}$ & 1 if and only if vehicle $k\in K$ uses time window $p\in P_i$ for visit $i\in N$ \\
$a^k_i\in \mathbb{R}$     & Time at which service starts at visit $i\in N$ \\
\hline
\end{tabular}
\end{center}
\caption{Mathematical notation\label{tab:notation}}
\end{table}

The following formulation is based on four group of decision variables: $x^k_{ij}\in \mathbb{B}$ indicating if vehicle $k\in K$ services visit $j\in V$ after visit $i\in V$. The variable $z^k_{ip}\in \mathbb{B}$ is $1$ if and only if vehicle $k$ services visit $i$ during time-window $p\in P_i$. The start of service and waiting times of the vehicles are decided by the variables $a^k_i\in \mathbb{R}$ and $w^k_i\in \mathbb{R}$ respectively.

\begin{align}
&\text{\bf minimise}\nonumber\\
&\sum_{k\in K}{ \sum_{i,j\in V}{ c^k_{ij}x^k_{ij} }  }+B\sum_{i\in N}{\left(s_{i}+\sum_{k\in K}{w_i^k}\right)}+
    \sum_{k\in K}{ c^k \left(\sum_{j\in N}{ x^k_{d_{start}j}} \right)}\label{eq:obj}\\
&\text{\bf subject to}\nonumber\\
&   \sum_{k\in K}{\sum_{j \in V}{x^k_{ij}} } = 1  \quad  \forall i\in N \label{eq:coverage}\\
&   \sum_{j \in V} x_{ji}^k = \sum_{j \in V} x_{ij}^k \quad  \forall i\in N, k\in K \label{eq:flow} \\
&   \sum_{j \in V}{x^k_{d_{start}j}} \leq 1  \quad  \forall k\in K  \label{eq:leave}\\
&   \sum_{j \in V}{x^k_{d_{end}j} } + \sum_{i \in V}{x^k_{i d_{start}} } = 0  \quad  \forall k\in K \label{eq:reenter}\\
&   \sum_{i \in V} \sum_{j\in N} q_j x_{ij}^k \leq Q_k \quad  \forall k\in K  \label{eq:capacity}\\
&   \sum_{j \in V} x_{ij}^k = \sum_{p \in P_i} z_{ip}^k \quad  \forall i\in V, k\in K  \label{eq:tw_choose}\\
& a^k_j \geq a^k_i + t_{ij} + s_i - M(1-x^k_{ij}) \quad  \forall k\in K, i,j\in V  \label{eq:service_time}\\
& a^k_i \geq l^p_i - M(1-z^k_{i,p}) \quad  \forall k\in K, i\in V, p\in P_i \label{eq:tw1}\\
& a^k_i \leq u^p_i + M(1-z^k_{i,p}) \quad  \forall k\in K, i\in V, p\in P_i \label{eq:tw2}\\
& w^k_j \geq a^k_j - a^k_i - t_{ij} - s_i - M(1-x^k_{ij}) \quad  \forall k\in K, i\in V ,j\in N \label{eq:waiting}\\
& T^k \geq a^k_{d_{end}} - a^k_{d_{start}} \quad  \forall k\in K  \label{eq:duration}\\
& w^k_i, a^k_i \geq 0 \quad  \forall i \in N, k \in K \label{eq:domainR}\\
& x^k_{ij}, z^k_{i,p} \in \{0,1\} \quad  \forall i,j \in N, k \in K, p\in P_i \label{eq:domainB}
\end{align}

The objective function~\eqref{eq:obj} minimises three terms: the total cost of the route, the total service and waiting time and the cost of vehicles used. The inclusion of the service and waiting time into the objective is dictated by the parameter $B$, which is $1$ if we want to allow this minimisation and $0$ otherwise. Equations~\eqref{eq:coverage} are the coverage constraints, while \eqref{eq:flow} are the flow balance constraints. We ensure that each vehicle leaves its start depot at most once~\eqref{eq:leave} and that it cannot return to it or leave the end depot~\eqref{eq:reenter}. The capacity of a vehicle is restricted by~\eqref{eq:capacity}. The modelling of the multiple time-windows begins with Constraints~\eqref{eq:tw_choose}. This constraint ensures that if a vehicle $k$ services visit $i$, it must choose one of the available time-windows. The selection will then be indicated by the variable $z^k_{ip}$. The time when the service of visit $i$ is started ($a^k_i$) is defined by~\eqref{eq:service_time}. The service start time is then constrained to be within the selected time-window with Constraints~\eqref{eq:tw1} and \eqref{eq:tw2}. The waiting time is calculated by~\eqref{eq:waiting} and the route duration is restricted by~\eqref{eq:duration}. Finally,~\eqref{eq:domainR} and \eqref{eq:domainB} define the domain of the variables.

The model is descriptive only as attempts to obtain bounds or solutions were unsuccessful.

\section{An \ac{ALNS} for the \ac{VRPMTW}}
\label{sec:ALNS}
Typically ALNS uses a solution representation consisting of a list of visits, (earliest/latest) arrival times at each visit, and (forward/backward) slack describing how much forward/backward in time each visit can be moved. The information added to the list allows fast evaluation of insertion costs/feasibility, and can trivially be extended to \ac{VRPMTW} by adding information about the time window chosen for each visit upon insertion.

We propose augmenting the list of visits with sets of labels for each visit, allowing efficient evaluation of the cheapest insertion for \emph{any} allocation of time windows. Crucially, infeasible insertions can be detected in constant time.

We test our solution representation on a version of the \ac{ALNS} framework from \citep{ALNS2006} adapted to the \ac{VRPMTW}. \ac{ALNS} is an extension of the Large Neighbourhood Search from \citep{Shaw1998}, where, using a destroy operator, a portion of the current solution is destructed. A repair heuristic is then applied to re-build the solution. \ac{ALNS} extends this search procedure by allowing multiple destroy and repair heuristics. At each iteration, it selects a destroy and a repair method using an adaptive engine. It applies the destroy method to remove visits, and applies the repair method to reinsert them. If the cost of the solution decreases, it is accepted. If it is not, the acceptance criterion from simulated annealing is applied, which results in a decreasing acceptance of worsening solutions as the search progresses, and makes acceptance increasingly unlikely as the cost increase grows. Afterwards, the temperature is adjusted using a geometric cooling schedule, and the likelihood of re-selecting the destroy and repair methods is updated based on the outcome of the destroy-repair operation. The procedure terminates once it reaches a time limit. We first present an overview of the algorithm, followed by a description of its components.
\begin{description}
\item[Step 1 - Initial Solution)]{A greedy insertion heuristic is used to find a first feasible solution to the problem. We randomly iterate over visits. Once we find a route where the visit can be inserted we apply a best insertion. Should no route be feasible a new one is opened. The algorithm continues until all visits are assigned.}

\item[Step 2 - Online temperature tuning] {Since the cooling schedule of the simulated annealing acceptance criterion can have a large impact on the search, we implemented a form of online instance specific automated tuning. To do so, we execute the destroy and repair part of our \ac{ALNS} for a limited number of iterations, allowing only improving solutions. During the search, the collected statistics are used to define an instance specific cooling schedule. More details about this step can be found in Section~\ref{sec:onlineTuning}.}

\item[Step 3 - Route minimisation] {In order to minimise the number of routes used, we artificially increase the cost of each route, forbid adding routes and add a cost for having unassigned visits. Starting from the initial greedy solution, we unassign all the visits of the smallest route. We then look for a new solution, using the destroy and repair part of the \ac{ALNS}, where the number of allowed routes is reduced by one. This process iterates until a time limit is reached. The best feasible solution that minimises the number of routes is then returned. A detailed description of the route minimisation step can be found in Section~\ref{sec:routeMinimisation}.}

\item[Step 4 - Solution optimisation] {In the final step, the solution from Step 3 is used as starting point. Here, the route costs are returned to their original values and unassigned visits are disallowed. The destroy and repair part of the \ac{ALNS} algorithm is run for the reminder of the time limit. }
\end{description}

Step 2-4 all utilise the same destroy and repair routines, differing only in the objective function, strategy for saving the best found solution and a strategy for when to remove entire routes routes (during route minimisation).

\subsection{Destroy methods}
We apply seven different destroy methods. For all of them, the number of visits to remove from the solution is chosen randomly from the interval of $[10,40]$. Each destroy method is then tasked with removing the visits from the solution. ALNS is a well known method for solving \ac{VRPTW}s, therefore a number of efficient destroy and repair methods have already been studied \citep{ALNS2006, ALNS2}. We apply them to \ac{VRPMTW}.
\begin{description}
\item [Random] from \cite{ALNS2006} selects the visits to remove randomly.
\item [Cluster $y$] for $y\in\{1,2,4\}$ are three destroy methods that select visits to remove randomly, but also remove $y$ successors to each visit. The method is based on ideas from the VeRoLog conference, but no publication of the method could be found.
\item [Geometric] adapted from \citep{ALNS2}, selects a random seed visit. It then repeatedly selects a removed visit $v_{select}$, sorts the visits remaining in the solution in ascending order by distance to $v_{select}$, and selects the visit that lies on the $r^4*n_{rem}$'th position of the ordered list, where $r$ is a random number in the interval $[0,1]$ and $n_{rem}$ is the number of visits remaining in the solution.
\item [Time] adapted from \citep{ALNS2}, selects visits like the geometric destroy, but sorts remaining visits in ascending order according to the difference between each visits start time and the start time of the randomly selected visit scheduled for removal $v_{select}$ in the current solution.
\item [Solution history] adapted from \citep{ALNS2}, is based on an assumption that if edges are not typically found in good solutions, a better alternative is more likely to exist, and removing them from a solution is more likely to yield improvements. It keeps track of the cost of the best solution observed, that contains any edge $(i,j)$. It scores each visit using sum of the score of the two edges connected to it in the solution. All visits are sorted in descending order according to their score, and the visits are removed choosing the $r^4 \cdot n_{rem}$'th where $r$ is a random number in the interval $[0,1]$.

Concretely for a visit $i$ in a solution, connected with edges $(pre(i),i)$ and $(i,suc(i))$ in the incumbent solution: If the best solution observed using the edge $(pre(i),i)$ cost 5000 and the best solution using edge $(i,suc(i))$ cost 4500, the score of the visit is 9500.
\end{description}

\subsection{Repair methods}
After each destroy method has been employed, one repair method from \citep{ALNS2006, ALNS2} is chosen and applied to all unassigned visits.

\begin{description}
\item [Regret-2] was presented in \citep{ALNS2006} who adapted methods from \citep{regret}. We evaluate the cost of inserting each unassigned visit in each route, and rank the visits in decreasing order of the cost of their second best insertion. The visit with the highest cost of insertion into its second best route, is inserted into its best route.
\item [Regret-2 randomised] is a randomised version of Regret-2. It ranks all unassigned visits in increasing order of the cost of their insertion into the second best route, but each cost is increased with a percentage chosen uniformly randomly in [0\%,50\%].
\end{description}

\subsubsection{Caching insertion costs}
Each time an insertion has been made, all unassigned visits have been evaluated for insertion in all routes. After the insertion, all previously calculated insertion costs are valid except those pertaining to the route wherein the visit was inserted. We employ this standard caching approach, where the insertion cost of each visit into each route is cached lazily. The cache is invalidated for all insertions pertaining to a route, each time a visit is inserted into the route.

\subsection{\texorpdfstring{$\Delta$}{Delta}-evaluation and feasibility}
ALNS, as well as most common heuristics, needs an efficient calculation of the cost of inserting a visit into a route in a given position. This is typically achieved using $\Delta$-evaluation (also known as incremental updates) \cite[pp48]{Hoos:2004:SLS:983505}, where only a (small) part of the full objective needs to be recalculated.
For a VRPTW this can be accomplished in constant time per position evaluated by caching enough information to establish an earliest possible start time, and a latest possible start time for each visit. The cost of maintaining these are typically $O(n)$ per update to the route, where $n$ is the size of the route. These methods can be generalised for both variants of the VRPMTW. Note that we include the start and the end depot in the routes and associate time windows to them.

\subsubsection{Without time minimisation}
\label{sec:delta_no_time}
In the literature, the VRPTW and VRPMTW variants without time minimisation are the most common. The relative ease in evaluating feasibility and cost of insertions has been observed by others, and makes this version of the problem a good starting point to understand the complexities of VRPMTW. We easily adapt forward slack by \citep{Savelsbergh} to work on \ac{VRPMTW} without time minimisation.

Let $pre(i)$ and $suc(i)$ define the predecessor and the successor respectively of the visit $i$ in a given solution. 

Let $es_i$ and $ls_i$ for all $i \in V$ be the earliest and latest start time of visit $i$ respectively. We have that

\begin{equation}\label{eq:es}
es_i = \min_{p\in P_i | u^p_i \geq es_{pre(i)}+t_{pre(i)i}} max(es_{pre(i)}+t_{pre(i)i},l^p_i) 
\end{equation}
with a base case of $es_i = \min_{p\in P_i} l^p_i$ if $i$ has no predecessors.
\begin{equation}\label{eq:ls}
ls_i = \max_{p\in P_i | l^p_i \leq ls_{suc(i)}-t_{isuc(i)}} min(ls_{suc(i)}-t_{isuc(i)},u^p_i) 
\end{equation}
with a base case of $ls_i = \max_{p\in P_i} u^p_i$ if $i$ has no successors.

Equation~\eqref{eq:es} and~\eqref{eq:ls} can be evaluated in $O(|P_i|)$, making the whole update procedure take $O(n|P_{max}|)$ where $P_{max}= \max_{i\in V} |P_i|$.

Insertion of $j$ between $i$ and $suc(i)$ in the incumbent solution using time window $p\in P_j$ in a route is feasible if and only if:
\begin{equation}\label{eq:insfeas}
es_i+t_{ij} \leq u^p_j \land ls_{suc(i)}-t_{jsuc(i)}\geq l^p_j
\end{equation}

and the cost of the insertion is given by:

\begin{equation}\label{eq:cost}
cost=c_{ij}+c_{jsuc(i)}-c_{isuc(i)}
\end{equation}

Note that Equation~\eqref{eq:insfeas} is verifiable in $O(|P_j|)$, and Equation~\eqref{eq:cost} is constant time. For the problems treated in this paper, the number of time-windows is limited, and further optimisations are not considered.

\subsubsection{With time minimisation}
When time minimisation is included in the objective, the calculations from Section~\ref{sec:delta_no_time} are insufficient to establish the cost of an insertion. Two distinct methods can be applied when treating this problem:
\begin{enumerate}
    \item When inserting a visit, keep track of and fix the time-window used for insertion.
    \item Fix only the ordering of the visits, and evaluate insertions for all possible allocations of time-windows to visits.
\end{enumerate}
The former suffers from routes being very hard to create good initial solutions for, as decisions regarding a visit expected position in the route is taken on insertion, including for the first inserted visit. Furthermore, difficulties with undoing the sub-optimal routes were observed during local search or destroy repair iterations, due to larger parts of the solution having to be changed in order to effect the change.

The second method suffers from the increased complexity in evaluating feasibility and cost of insertions. These challenges however, are not insurmountable.

When evaluating the insertion of a visit on a given position in a route, feasibility can be evaluated as in the case without time minimisation. To evaluate the resulting minimum route duration, we apply a labelling algorithm which is a version of dynamic programming. It works by keeping track of any promising (non dominated) ways the vehicle could arrive at the current visit. This includes when the vehicle should start, implicitly what time windows it should use on the way, and how much later the vehicle could start and still hit the time windows. Similar labels are kept describing any possible path from each visit, along the rest of the route and back to the depot: For each visit we maintain two sets of labels $L_i^f$ and $L_i^b$, namely the labels describing forwards and backwards paths from each visit $i\in V$. $L_i^b$ consists of labels containing the earliest start time $es_i$, the cost-less backwards slack $bs_i$, and resulting route start time $st_i$. $L_i^f$ consists of labels containing the latest start time $ls_i$, the cost-less forward slack $fs_i$, and the resulting route end time $et_i$.

Given two consecutive visits in a route $i$ and $j$, any backward labels $(es_i,bs_i,st_i)\in L_i^f$ and any forward label $(ls_j,fs_j,et_j)\in L_j^b$. The earliest start time for the visit we wish to insert ($o$) using the $p$'th time window $[l^p_o,u^p_o]$ is $es_o = max(es_i+t_{io},l^p_o)$ and the latest start time is $ls_o = min(ls_j-t_{oj}-s_o,u^p_o)$. The insertion is feasible if and only if $ls_o>es_o$. The resulting duration of the route is $et_j-st_i-min(ls_o-es_o, bs_i+fs_j)$. The term $min(ls_o-es_o, bs_i+fs_j)$ is due to compacting of the route by time shifting the forward and backward paths.

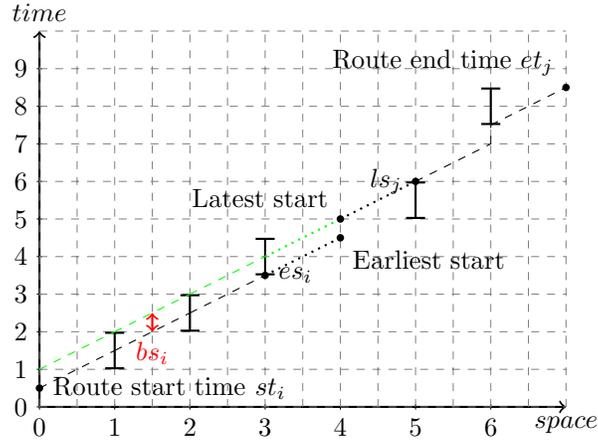
\begin{figure}[ht]
    \centering

\begin{tikzpicture}[scale=1,dot/.style 2 args={circle,inner sep=1pt,fill,label={#2:#1},name=#1},
  extended line/.style={shorten >=-#1,shorten <=-#1,draw=Cerulean},
  extended line/.default=1cm]
    \coordinate (A) at (7,0);
    \coordinate (B) at (0,5);

    \draw[thick,->] (0,0) -- (A) node[anchor=north] {$space$};
    \draw[thick,->] (0,0) -- (B) node[anchor=south] {$time$};
    
    \foreach \x in {0,1,2,3,4,5,6}
        \draw (\x cm,1pt) -- (\x cm,-1pt) node[anchor=north] {$\x$};
    
    \foreach \y in {0,1,2,3,4,5,6,7,8,9}
        \draw (1pt, \y*0.5 cm) -- (-1pt,\y*0.5 cm) node[anchor=east] {$\y$};
    
    \draw[step=0.5cm,gray, dashed ,very thin] (0,0) grid (7,5);

    \draw[thick,|-|] (1,0.5) coordinate (a_1) -- (1,1) coordinate (b_1);
    \draw[thick,|-|] (2,1) coordinate (a_2) -- (2,1.5) coordinate (b_2);
    \draw[thick,|-|] (3,1.75) coordinate (a_3) -- (3,2.25) coordinate (b_3);
    
    \draw[thick,|-|] (5,2.5) coordinate (a_5) -- (5,3) coordinate (b_5);
    \draw[thick,|-|] (6,3.75) coordinate (a_6) -- (6,4.25) coordinate (b_6);

    \draw[dashed,green] (0,0.5) coordinate (a_1) -- (1,1) coordinate (a_2);    
    \draw[dashed,green] (1,1) coordinate (a_1) -- (2,1.5) coordinate (a_2);
    \draw[dashed,green] (2,1.5) coordinate (a_1) -- (3,2) coordinate (a_2);
    \draw[dotted,thick,green] (3,2) coordinate (a_1) -- (4,2.5) coordinate (a_2);
    
    \draw[dashed] (0,0.25) coordinate (a_1) -- (1,0.75) coordinate (a_2);    
    \draw[dashed] (1,0.75) coordinate (a_1) -- (2,1.25) coordinate (a_2);
    \draw[dashed] (2,1.25) coordinate (a_1) -- (3,1.75) coordinate (a_2);
    
    \draw[dotted,thick] (3,1.75) coordinate (a_1) -- (4,2.25) coordinate (a_2);
    \draw[dotted,thick] (4,2.5) coordinate (a_1) -- (5,3) coordinate (a_2);
    
    \draw[dashed] (5,3) coordinate (a_1) -- (6,3.5) coordinate (a_2);
    \draw[dashed] (6,3.5) coordinate (a_1) -- (6,3.75) coordinate (a_2);
    \draw[dashed] (6,3.75) coordinate (a_1) -- (7,4.25) coordinate (a_2);
    
    \node [dot={Route start time $st_i$}{right}] at (0,0.25) {};
    \node [dot={Route end time $et_j$}{above left}] at (7,4.25) {};

    \node [dot={$es_i$}{right}] at (3,1.75) {};
    \node [dot={$ls_j$}{left}] at (5,3) {};
    
    \node [dot={Earliest start}{below right}] at (4,2.25) {};
    \node [dot={Latest start}{above left}] at (4,2.5) {};
    
    \draw[thick, red,<->] (1.5,1.25) coordinate (a_6) -- (1.5,1) coordinate (b_6) node[anchor=north] {$bs_i$};

\end{tikzpicture}

    \caption{Forward ($L^f_{j}=(es_j,bs_j,st_j)$) and backward ($L^b_{i}=(ls_i,fs_i,et_i)$) paths/labels for evaluating insertion at space location indexed by 4.}
    \label{fig:twPaths}
\end{figure}

An example of inserting visit $o$ between visits $i$ and $j$ using the labels is illustrated in figure \ref{fig:twPaths}. The backwards path (black, dashed) from space index 0 to space index 3 indicates the driving time to get to the visit at space index 3 as early as possible, without starting the route earlier than necessary ($st_i$). The path can be shifted forward in time (up to the green and dashed line) by $bs_i$ with a resulting change in start time up to $st_i+bs_i$. Further delays are possible, but the resulting start time will be capped at $st_i+bs_i$. The forward path from space index 5 to 7 (black, dashed) is the latest possible path to arrive at the earliest possible end time $et_j$. If this is shifted earlier, the resulting end time will remain unchanged, as the path is already constrained by the time window at space index 6, an thus contains no slack ($fs_j=0$).

\subsubsection{Label generation and elimination}
Let the start depots backwards label set $L^b_{d_{start}}$ consist of the label $(0,\infty,0)$ and the end depots forward label set $L^f_{d_{end}}$ consist of the label $(\infty,\infty,\infty)$. With these base cases we define the remaining sets recursively.

\begin{description}
\item[Backward labels:] Let $j=pre(i)$ and $(es_j,bs_j,st_j) \in L^b_{j}$ be a backward label of $j$. Using the $p$'th time window $[l^p_i,u^p_i]$ a new label $(es_i, bs_i, st_i)$ can be generated as $es_i = max(es_j+t_{ji}+s_{j}, l_i^p)$, $bs_i=max(0,bs_j-max(0,l^p_i-es_j-t_{ji}-s_j))$ and $st_i=st_j+bs_j-bs_i$.
\item[Forward labels:] Let $j=suc(i)$ and $(ls_j,fs_j,et_j) \in L^f_{j}$ be a forward label of $j$. Using the $p$'th time window $[l^p_i,u^p_i]$ a new label $(ls_i, fs_i, et_i)$ can be generated as $ls_i = min(ls_j-t_{ij}-s_{i}, u_i^p)$,
$fs_i=max(0,fs_j-max(0,ls_j+t_{ij}+s_i-u^p_i)$ and $et_i=et_j-max(0,ls_j+t_{ij}+s_i+u^p_i)$.
\end{description}

Given backward labels $(es_i,bs_i,st_i)$ and $(es_j,bs_j,st_j)$, we say that the former dominates the latter if and only if $st_i+bs_i \geq st_j+bs_j \wedge es_i \leq es_j$.

Given forward labels  $(ls_i,fs_i,et_i)$ and $(ls_j,fs_j,et_j)$, we say that the former dominates the latter if and only if $et_i-fs_i \leq et_j-fs_j \wedge ls_i \geq ls_j$.

\begin{algorithm}[h!] 
\KwData{A list of backward labels $L^b_{i}$ for a visit $i$, a set of forward labels $L^f_{j}$ for its successor $j$ and a visit $o$ to insert. Additionally the old duration of the route $oldDur$ is needed.}
\KwResult{The cheapest insertion delta cost of $o$ between $i$ and  $j$}
$bestCost=\infty$\;
\For{$(es_i,bs_i,st_i)\in L^b_i$ in increasing order of $es_i$}{ \label{alg:iterBL}
    \For{$(ls_j,fs_j,et_j)\in L^f_j$ in decreasing order of $ls_j$}{ \label{alg:iterFL}
        \If{$es_i+t_{io}\leq ls_j-t_{oj}$}{\label{alg:compat}
            \For{Each time window $[l^p_o,u^p_o] \in P_o$ in increasing order of time}{\label{alg:forTW}
                \If{$es_i+t_{io}\leq u^p_o$ and $ls_j-t_{oj}\geq l^p_o$}{ \label{alg:cond}
                    Compute label $(es_o,bs_o,st_o)$ by expanding $(es_i,bs_i,st_i)$ using $[l^p_o,u^p_o]$\; \label{alg:exp1}
                    Compute label $(ls_o,fs_o,et_o)$ by expanding $(ls_j,fs_j,et_j)$  using $[l^p_o,u^p_o]$\;\label{alg:exp2}
                    \If{$es_o \leq ls_o$}{  \label{alg:cond2}
                    $newDur=et_o-st_o-min(ls_o-es_o, bs_i+fs_i)$\; \label{alg:newCalc1}
                    $\Delta travelCost=t_{io}+t_{oj}-t_ij{}$\;
                    $\Delta cost= oldDur-newDur+travelCost$\; \label{alg:newCalc2}
                    \If{$\Delta cost<bestCost$}{
                        $bestCost=\Delta cost$
                    }
                    }
                }
            }
        }
    }
}
\Return{$bestCost$}\;
\caption{Algorithm finding the cheapest insertion cost for inserting a visit $o$ in a given index in a route.}
\label{alg:checking}
\end{algorithm}

The cheapest insertion of visit $o$ between visits $i$ and $j$ can the be calculated using Algorithm \ref{alg:checking}. It iterates through the backward labels for the previous visit $i$ (line \ref{alg:iterBL}) and the forward labels for the next visit $j$ (line \ref{alg:iterFL}), checking for compatibility ($es_i+t_{io}>ls_j-t_{oj}$) (line \ref{alg:compat}). If the labels are compatible in the sense that driving time seems feasible without time windows, they are expanded for each matching time window (line \ref{alg:forTW}) to visit $o$ as if it was already inserted (line \ref{alg:exp1}-\ref{alg:exp2}). The new forward and backward labels for $o$ are (if $es_o \leq ls_o$) then used to compute the route duration of the resulting route, and the delta cost of the insertion (line \ref{alg:newCalc1}-\ref{alg:newCalc2}). By carefully considering the conditions in line \ref{alg:compat}, \ref{alg:cond} and \ref{alg:cond2}, iterations through either of the sets of labels can often be elided. This is crucial to obtain the performance observed in the paper.

\subsection{Adaptive engine}
ALNS uses an adaptive engine from \citep{ALNS2006} to select between destroy methods $D$ and repair methods $R$. The probability of selecting a specific destroy method $de$ is $\frac{weight_{de}}{\sum_{des\in D} weight_{des}}$, and the probability of selecting a repair method $re$ is $\frac{weight_{re}}{\sum_{res\in R} weight_{res}}$.

Each time a repair or destroy method is chosen, its weight is updated using a score and a decay factor $decay$. The score $sc$ depends on the result of the destruction and reconstruction, and is set to 1, $score_a$, $score_i$ or $score_b$ depending on if the new solution is rejected, accepted, improving or a new best respectively. The scores and the decay factors are parameters to be set during tuning.
The new weight is calculated as $weight_{new}=weight_{old}\cdot decay +  sc \cdot (1-decay)$.

The weight selection ensures that no destroy and repair method are ever entirely excluded as its weight can never be reduced below 1.

\subsection{Online temperature tuning}
\label{sec:onlineTuning}
We use an automatic tuning to control the simulated annealing part of the ALNS. The online temperature  is parametrised by a number of iterations $n_{init}$ to run only accepting improving solutions, as well as $\Delta cost_{init}$ and $\Delta cost_{end}$ the desired worsening of the solution to be accepted with a 50\% probability during the start and the end of the search respectively. The initial iterations establishes an approximation for the $cost$ of the solutions to be used in setting the temperature.

Recall that in simulated annealing the acceptance criteria for solutions worse than the current is controlled by $r<e^{\frac{\Delta cost}{cost}\cdot \tau}$, where $r\in[0,1]$ is a random number from a uniform distribution, and $\tau$ is the current temperature.

After the first $n_{iterTT}$ iterations, the start and the end temperature is determined.
\begin{equation}
\label{eq:starttemp}
\tau_{start}=\frac{log(0.5) \cdot cost_{init}}{\Delta cost_{init}}, \tau_{end}=\frac{log(0.5) \cdot cost_{end}}{\Delta cost_{end}}
\end{equation}
They are computed using \eqref{eq:starttemp} with the given target $\Delta cost_{init}$ and $\Delta cost_{end}$. A geometric cooling scheme is used to adjust the temperature $\tau$ during the search.

\subsection{Route minimisation}
\label{sec:routeMinimisation}

The strategy for minimising the number of routes, is based on the concept presented in \citep{ALNS2006}. An initial solution is created using the greedy construction heuristics in Step 1. After the temperature tuning (Step 2), the remaining optimisation is performed using the destroy and repair part of the \ac{ALNS} described in Section \ref{sec:ALNS}, with the following modifications.

The optimisation is divided into two phases: A route minimisation phase (Step 3) and an optimisation phase (Step 4). During the route minimisation phase, unassigned visits are penalised with a cost of $p_u=10000$ and routes are penalised with a cost of $p_r=1000000$. Thus any partial solution with fewer routes will be better cost wise if less than $\frac{p_r}{p_u}$ visits are unassigned. \ac{ALNS} is prevented from storing a partial solution as the best solution $s_{best}$, unless all visits have been inserted into the solution. Each time a solution with all visits assigned is obtained, the route with the fewest visits is destroyed.

\begin{figure}[ht]

    \centering

\begin{tikzpicture}[scale=1,dot/.style 2 args={circle,inner sep=1pt,fill,label={#2:#1},name=#1},
  extended line/.style={shorten >=-#1,shorten <=-#1,draw=Cerulean},
  extended line/.default=1cm]
    \coordinate (A) at (5,0);
    \coordinate (B) at (0,3);

    \draw[thick,->] (0,0) -- (A) node[anchor=west] {Percent done};
    \draw[thick,->] (0,0) -- (B) node[anchor=south] {Unassigned accepted};
    
    \foreach \x in {0,10,20,30,40,50,60,70,80,90,100}
        \draw (\x*0.05 cm,1pt) -- (\x*0.05 cm,-1pt) node[anchor=north] {$\x$};
    
    \foreach \y in {0,1,2,3,4,5}
        \draw (1pt, \y*0.5 cm) -- (-1pt,\y*0.5 cm) node[anchor=east] {$\y$};
    
    \draw[step=0.5cm,gray, dashed ,very thin] (0,0) grid (5,3);

    \draw[thick,dashed,->] (0.5,2.5) coordinate (a_1) -- (0.5,3) coordinate (a_2);

    \draw[thick] (0.5,2.5) coordinate (a_1) -- (2,2.5) coordinate (a_2);
    \draw[thick] (2,2.5) coordinate (a_1) -- (2,2) coordinate (a_2);

    \draw[thick] (2,2) coordinate (a_1) -- (2.5,2) coordinate (a_2);
    \draw[thick] (2.5,2) coordinate (a_1) -- (2.5,1.5) coordinate (a_2);

    \draw[thick] (2.5,1.5) coordinate (a_1) -- (3,1.5) coordinate (a_2);
    \draw[thick] (3,1.5) coordinate (a_1) -- (3,1) coordinate (a_2);

    \draw[thick] (3,1) coordinate (a_1) -- (3.5,1) coordinate (a_2);
    \draw[thick] (3.5,1) coordinate (a_1) -- (3.5,0.5) coordinate (a_2);
    
    \draw[thick] (3.5,0.5) coordinate (a_1) -- (4,0.5) coordinate (a_2);
    \draw[thick] (4,0.5) coordinate (a_1) -- (4,0) coordinate (a_2);

\end{tikzpicture}

    \caption{The local search will stay in route minimisation mode as long as the number of unassigned visits remain below the plotted line.}
    \label{fig:minim}
\end{figure}
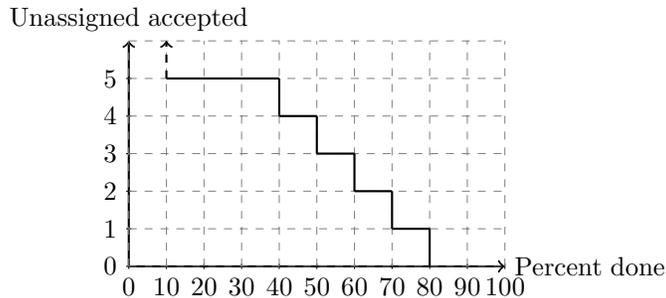

The shift from the route minimisation phase (Step 3) to the optimisation phase (Step 4) is triggered when the number of visits being unassigned is too high for the stage of the optimisation as depicted in Figure \ref{fig:minim}. This corresponds to one of the following conditions are fulfilled: 1) Ten percent of the allotted time has passed, and no visits are unassigned in the current solution. 2) Ten percent of the time has passed, and more than 5 visits remain unassigned. (See parameter tuning in Section \ref{sec:params}) 3) Ten percent of the time has passed, more than $percent/10-2$ visits remain unassigned, where $percent$ indicates the percentage of computation time is left. 
In the two latter cases, \ac{ALNS} failed to make the partial solution feasible with the given number of routes, and then algorithm reverts to the best found solution $s_{best}$ that is guaranteed to be feasible.

\section{Configuring the algorithm}
\label{sec:config}
The \ac{ALNS} just described, is composed of a set of components, each of which may have parameters to set, or may be omitted entirely. In this section we provide our rationalisation for the chosen parameters, and justify each components inclusion if it can be omitted.

Optimally each possible combination of parameters should be explored, along with each possible alternative implementation. Due to the sheer impossibility of this, a reduction in the search space is necessary. Where possible, we adopt values and methods known to work on VRPTW, and analyse deviations to increase the odds of at least a local optimum in the configuration space.

Initially we fixed the interval, from which the  number of visits to release is chosen. Intervals start and end values based on ($[10,40]$ from \citep{ALNS2}) were investigated in increments of 5, but little change was observed around the chosen values, with a sharp degeneration in solution quality if the values increased by more than 15.

Next the strategy transitioning from route minimisation (Step 3) to solution optimisation (Step 4). Initial trials indicated a threshold of 5 unassigned customers after 10\% of the search time had been expended, could be resolved in some cases. It was further observed that most subsequent improvements (Step 4) occur within the first 20\% of the optimisation time, and thus this amount should be reserved for the solution optimisation (Step 4). A linear decline from 40\% was chosen provisionally. Subsequent attempts to tune these parameters were unsuccessful due to the interdependence of the variables, and the big cost difference between solutions using varying number of vehicles being triggered both randomly and in a vast minority of instances.


\subsection{Parameter tuning}
\label{sec:params}

The algorithm is associated with a set of parameters to be tuned. To accomplish the tuning, we modify the instances by randomly reassigning the sets of time windows to the visits, thus creating a training set.

 We have limited ourselves to setting two sets of parameters for all instances. One for the problem with time minimisation, and a set of parameters for the problem without time  minimisation.

 \ac{ALNS} selects its repair and destroy methods using scores for previous performance. Recall that we apply the traditional selection procedures from \citep{ALNS2006}, and score $score_{b}$,  $score_{i}$ and  $score_{a}$ points for a new best, an improvement to the current solution and an accepted solution respectively. The decay of the score $decay\in [0,1]$ determines the memory of the adaptive engine. Eg. $0.9$ meaning 10\% of the score is based on the methods last score, and 90\% on the previous score.
 
Parameters are tuned using the {\em irace} package from \citep{LopDubPerStuBir2016irace}. Irace works by sampling the domain of each parameter and running the resulting algorithm configurations against each other. Each tuning is given a budget of 15000 algorithm runs, and it will statistically eliminate configurations as they prove inferior. Amongst the survivors the best performing configurations are presented.

\subsubsection{Tuning without time minimisation}
An initial run of the tuning was done to establish a baseline for each parameter, and test any assumptions made regarding minimum and maximum values of each parameter. These domains are specified in first line of table \ref{tab:initNoMin}.

\begin{table}[H]
    \begin{tabular}{r|r|r|r|r|r|r}
    $n_{iterTT}$ & $\Delta cost_{init}$ & $\Delta cost_{end}$ & $score_{a}$ & $score_{i}$ & $score_{b}$ & $decay$ \\ \hline
[10,$10^4$]& [0,1000] & [0,100] & [1,20] & [1,20] & [1,20] & [0.7,0.99]\\ \hline \hline
5641 & 213.93 & 4.38 & 4 & 10 & 15 & 0.77\\ \hline
5470 & 15.50 & 2.75 & 4 & 10 & 14 & 0.74\\ \hline
4801 & 201.14 & 1.16 & 4 & 10 & 15 & 0.75\\ \hline
4405 & 812.83 & 1.07 & 7 & 16 & 19 & 0.84\\ \hline
5312 & 46.97 & 4.57 & 4 & 10 & 14 & 0.74\\
    \end{tabular}
    \caption{Initial tuning for the VRPMTW without time minimisation.}
    \label{tab:initNoMin}
\end{table}

The initial run resulted in the 5 configuration settings in table \ref{tab:initNoMin}. All values are well within the specified domains rather and none of them are on an extreme value. In order to further converge the parameter settings, already established values $n_{iterTT}$ and $score_{a}$ are fixed. This also serves to remove equivalence causing symmetries from the set of algorithm configurations. All remaining parameter domains are further tightened to the values in Table \ref{tab:endNoMin}.

\begin{table}[H]
    \begin{tabular}{r|r|r|r|r|r|r}
    $n_{iterTT}$ & $\Delta cost_{init}$ & $\Delta cost_{end}$ & $score_{a}$ & $score_{i}$ & $score_{b}$ & $decay$ \\ \hline
5000 & [10,800] & [0,5] & 4 & [10,20] & [10,20] & [0.7,0.99]\\ \hline \hline
5000 & 298.49 & 2.24 & 4 & 15 & 15 & 0.73\\ \hline
5000 & 223.80 & 1.78 & 4 & 16 & 20 & 0.73\\ \hline
5000 & 393.68 & 2.73 & 4 & 14 & 14 & 0.75\\ \hline
5000 & 332.21 & 2.58 & 4 & 15 & 18 & 0.76\\ \hline \hline 
5000 & 300.00 & 2.33 & 4 & 15 & 17 & 0.75\\
    \end{tabular}
    \caption{Final tuning for the VRPMTW without time minimisation. The final line is the chosen configuration.}
    \label{tab:endNoMin}
\end{table}

Based on the parameter setting from in the top part of Table \ref{tab:endNoMin} a the parameter tuning recommended the parameters in the final row of the table.

Comparing these values to other ALNS parameter settings, they are mainly characterised by an extreme emphasis on diversification. Both in terms of acceptance of worsening solutions by having a high $\Delta cost_{init}$, but also in a diverse selection of destroy/repair methods in having a low $decay$ factor, keeping more emphasis on recent results than historical.

\subsubsection{Tuning with time minimisation}

The initial tuning/check was given the domains from Table \ref{tab:initMin}, and the following configurations were obtained:

\begin{table}[H]
    \begin{tabular}{r|r|r|r|r|r|r}
    $n_{iterTT}$ & $\Delta cost_{init}$ & $\Delta cost_{end}$ & $score_{a}$ & $score_{i}$ & $score_{b}$ & $decay$ \\ \hline
[10,$10^4$]& [0,1000] & [0,100] & [1,20] & [1,20] & [1,20] & [0.7,0.99]\\ \hline \hline
1028.33 & 19.46 & 4.91 & 2 & 11 & 14 & 0.93\\ \hline
1594.53 & 28.27 & 3.60 & 3 & 9 & 12 & 0.90\\ \hline
3942.04 & 518.73 & 2.14 & 7 & 11 & 14 & 0.90\\ \hline
570.89 & 20.63 & 6.02 & 2 & 11 & 14 & 0.92\\
    \end{tabular}
    \caption{Initial tuning with time minimisation}
    \label{tab:initMin}
\end{table}

From Table \ref{tab:initMin} a set of parameter domains for the final tuning was established and presented in the first row of Table \ref{tab:finalTuneMin}.

\begin{table}[H]
    \begin{tabular}{r|r|r|r|r|r|r}
    $n_{iterTT}$ & $\Delta cost_{init}$ & $\Delta cost_{end}$ & $score_{a}$ & $score_{i}$ & $score_{b}$ & $decay$ \\ \hline
[500,5000] & [10,600] & [0,10] & [1,7] & [1,20] & [10,20] & [0.8,0.99]\\ \hline \hline
1999 & 44.65 & 4.95 & 2 & 9  & 17 & 0.83\\ \hline
1590 & 14.32 & 6.19 & 2 & 6  & 16 & 0.84\\ \hline
1219 & 11.16 & 5.31 & 3 & 11 & 12 & 0.84\\ \hline
1530 & 35.16 & 6.14 & 2 & 6  & 17 & 0.83\\ \hline \hline
1600 & 26.00 & 5.50 & 2 & 8  & 16 & 0.83\\
    \end{tabular}
    \caption{Final tuning with time minimisation}
    \label{tab:finalTuneMin}
\end{table}

The four candidate configurations presented in Table \ref{tab:finalTuneMin} were reduced to the final algorithm configuration in the final row. These values continue the trend of diversification observed for the VRPMTW without time minimisation, but is adjusted for the increased time usage per iteration by having fewer pre-iterations ($n_{iterTT}$), and much lower start temperature.

\subsubsection{Importance of tuning}
Note that some parameters such as $decay$ and $n_{iterTT}$ converge faster during the tuning, and thus can be assumed to have more impact on the final solution quality. Other parameters such as $\Delta cost_{init}$ seems to vary more. This lack of homogeneity between the configurations can be due to (a combination of) either a lack of tests to create statistical significance or lack of meaningful impact on the solution quality. In many applications tuning have shown to provide a great advantage, but our original test runs using the default configurations from Table \ref{tab:defaultConf} gave solutions with only marginally worse average solutions, and one more best-known solutions. We thus conclude that further tuning is unwarranted.
\begin{table}[H]
    \begin{tabular}{r|r|r|r|r|r|r}
    $n_{iterTT}$ & $\Delta cost_{init}$ & $\Delta cost_{end}$ & $score_{a}$ & $score_{i}$ & $score_{b}$ & $decay$ \\ \hline
1000 & 10.00 & 3.0 & 2 & 4  & 10 & 0.9 \\
    \end{tabular}
    \caption{Original parameter values}
    \label{tab:defaultConf}
\end{table}

\subsection{The value of the solver components}
To evaluate the importance of the solver components, we have omitted any component not strictly required. For implicit time windows (e.g. evaluating the new cost for the optimal allocation of time windows in the route) we implemented a modified solution representation that always inserts using the cheapest time window, and only changes time window upon reinsertion. When omitting online temperature tuning, we ran the algorithm with online temperature tuning initially, and fixed the average start and end temperatures selected. Each destroy and repair method can be omitted without complications.

\begin{figure}[htbp]
    \centering
    \begin{subfigure}[t]{0.4\textwidth}
        \includegraphics[width=\textwidth]{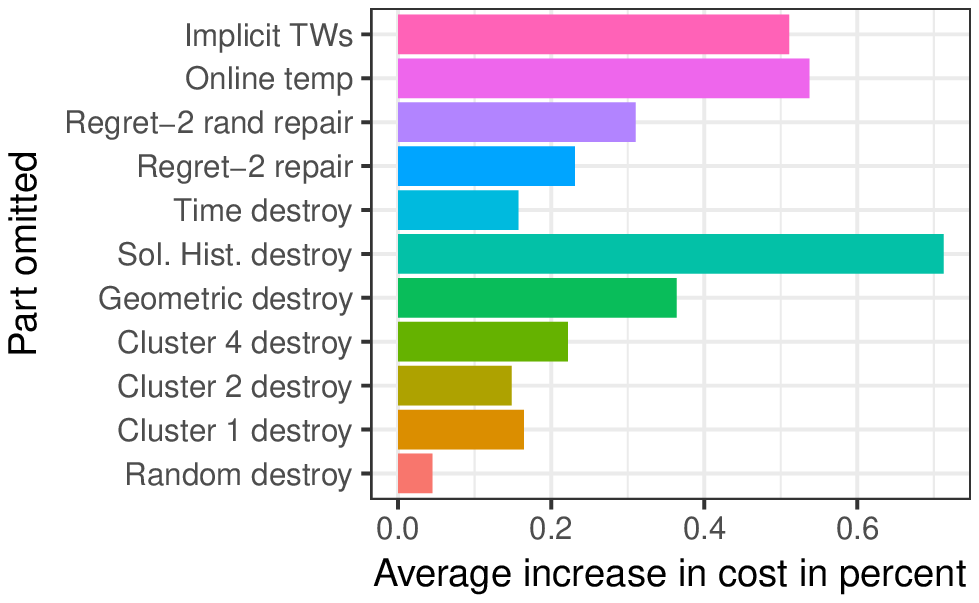}
        \caption{Increase in average cost}
        \label{fig:omit}
    \end{subfigure}
    \begin{subfigure}[t]{0.4\textwidth}
        \includegraphics[width=\textwidth]{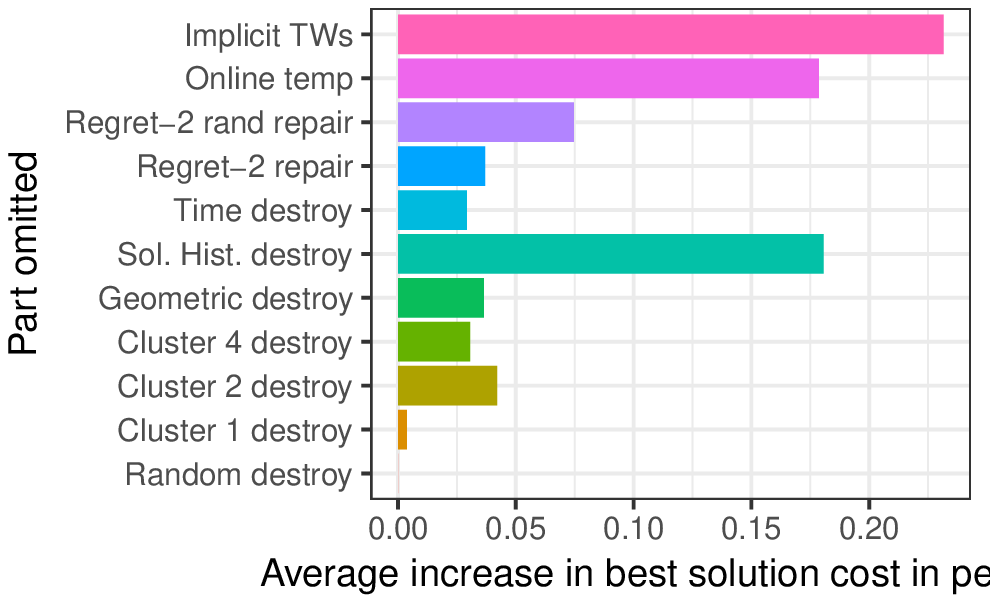}
        \caption{Increase in best solution cost}
        \label{fig:omitBest}
    \end{subfigure}
    \caption{Increase in solution cost when omitting components.}
\end{figure}

Figure \ref{fig:omit} indicates the increase in the average solution cost when running the time minimising algorithm for 60 seconds while omitting each of the tested components, while Figure \ref{fig:omitBest} shows the increase in the best solution cost for the same settings. Since in the instances most costs components such as service times, vehicle costs and minimum distances are of greater magnitude than travelled distances and waiting times, small changes are more significant than they appear.

The solution history (Sol. Hist.) destroy is the most significant component when looking at average costs in Figure \ref{fig:omit}, which may be due to expanded options for choosing time windows in order to place visits close, and thus a greater chance that close visits can be place close together. Online temperature tuning is artificially important as instances can be divided into cheap instances with costs less than 5000 and instances with costs in excess of 10000, each requiring different temperatures. Third most important is our implicit time windows. The sharp increase in iterations resulting from deactivating them is far from able to compensate for the loss of flexibility.


Figure \ref{fig:omitBest} indicates an increased importance of the implicit time window selection when looking for the best solutions rather than averages costs.

\section{Computational results}
\label{sec:results}
We test our adaptation of the \ac{ALNS} on the instances proposed by \cite{Belhaiza2014AWindows}. The instances are generated based on the traditional Solomon instances \citep{solomon} with a set of newly generated time windows. As a result the instances have 100 customers, and vehicle capacity of 200 (cm1xx, rcm1xx, rm1xx), 700 (cm2xx) or 1000 (rcm2xx,rm2xx)
All tests are performed on a Intel Core i7-4790K CPU @ 4.00GHz using a single thread. The results generated within a time limit of 60 seconds, which is slightly less than the average time used in \citep{Belhaiza2014AWindows}(they do not operate with a fixed time limit).

For the instances, the cost is equivalent to the time for traversing the edge between the visits $c^k_{ij}=t_{ij}$, and the cost of using a vehicle $c^k$ is equal to its capacity $Q_k$. We further restricted the problem to assume $a^k_{d^k_{end}}\leq D^k$, meaning the route duration constraint is further restricted to be a deadline for returning to the depot. Any solution to this restricted problem is also a solution to the original as $a^k_{d^k_{end}}\leq D^k \Rightarrow a^k_{d^k_{end}}-a^k_{d^k_{start}}\leq D^k$ if and only if $a^k_{d^k_{start}}>0$.

 \begin{table}[htbp]
 \centering
 \footnotesize
 \begin{tabular}{c c | r r | r r | r r } 
  &  & \multicolumn{2}{c|}{\cite{Belhaiza2014AWindows}}  & \multicolumn{4}{c}{ALNS}  \\
 \hline
 Instance & m & $B=0$  & & $B=0$  & 60 sec. & $B=0$  & 600 sec. \\
  &  & Best R. & Avg.  & Best R. & Avg. & Best R. & Avg. \\
 \hline
 
 rm101 & 10     & 2977.2 & 3005.0 & \bf{2968.8}   & \bf{2972.3} & \textbf{2968.8} & \textbf{2969.0}\\
 rm102 & 9      & 2759.4 & 2759.4 & \bf{2705.9}   & \bf{2721.4} & \textbf{2705.9} & \textbf{2716.7}\\
 rm103 & 9      & 2692.5 & 2710.5 & \bf{2680.4}   & \bf{2686.2} & \textbf{2680.4} & \textbf{2681.6}\\ 
 rm104 & 9      & 2696.6 & 2719.2 & \bf{2690.7}   & \bf{2691.8} & \textbf{2690.7} & \textbf{2691.5}\\
 rm105 & 9      & 2688.8 & 2711.0 & \bf{2683.7}   & \bf{2685.7} & \textbf{2683.7} & \textbf{2683.7}\\
 rm106 & 9      & \bf{2691.9} & \bf{2817.2} & 2700.9 & 2704.4   &          2700.9 & \textbf{2703.9}   \\
 rm107 & 9      & 2690.8 & 2714.7 & \bf{2685.1}   & \bf{2685.8} & \textbf{2685.1} & \textbf{2685.1}\\
 rm108 & 9      & 2729.1 & 2729.1 & \bf{2716.0}   & \bf{2723.9} & \textbf{2716.0} & \textbf{2717.4}\\

 \hline \hline

 rm201 & \bf{2} & 3711.4 & 3720.5 & \bf{2750.0} & \bf{2769.0} & \textbf{2745.6} & \textbf{2756.0}\\
 rm202 & 2      & 2698.1 & 2717.3 & \bf{2681.3} & \bf{2685.0} & \textbf{2681.3} & \textbf{2685.0}\\
 rm203 & 2      & 2686.1 & 2702.0 & \textit{\bf{2678.8}} & \bf{2682.2} & 2679.0 & \textbf{2679.7}\\ 
 rm204 & 2      & 2680.5 & 2691.2 & \bf{2676.1} & \bf{2677.3} & \textbf{2672.8} & \textbf{2673.1}\\
 rm205 & 2      & \bf{2671.0}&2688.2&\bf{2671.0}& \bf{2671.6} & \textbf{2671.0} & \textbf{2671.1}\\
 rm206 & 2      & 2686.3 & 2704.9 & \bf{2679.0} & \bf{2681.9} & \textbf{2679.0} & \textbf{2679.0}\\
 rm207 & 2      & 2678.2 & 2696.2 & \bf{2674.7} & \bf{2679.7} & \textbf{2674.7} & \textbf{2678.2}\\
 rm208 & 2      & \bf{2673.9} & 2690.7 & 2675.7 & \bf{2677.7} & \textbf{2673.9} & \textbf{2676.2}\\

 \hline \hline

 cm101 & 10     & \bf{3089.2} & \bf{3102.4} & 3151.4 & 3212.1 & \textbf{3073.0} & 3126.8\\
 cm102 & 12     & \bf{3426.9} & \bf{3426.9} & 3444.6 & 3488.9 & \textbf{3392.5} & 3463.0\\
 cm103 & \bf{11}& 3532.7 & 3572.7 & \bf{3465.3}& \bf{3526.2}  & \textbf{3434.6} & \textbf{3477.1}\\ 
 cm104 & \bf{13}& 4051.3 & 4058.0 & \bf{3892.9}& \bf{3922.3}  & \textbf{3859.6} & \textbf{3907.1}\\
 cm105 & \bf{10}& \bf{3060.6} & \bf{3077.3} & 3062.3 & 3093.9 & \textbf{3037.9} & \textbf{3051.7}\\
 cm106 & 10     & 2992.4 & 3020.2 & \bf{2982.2} & \bf{2990.3} & \textbf{2980.2} & \textbf{2983.3}\\
 cm107 & \bf{10}& 3256.5 & 3292.3 & \bf{3077.9}& \bf{3079.4}  & \textbf{3077.9} & \textbf{3077.9}\\
 cm108 & 10     & 2968.7 & 2973.1 & \bf{2965.4}  &\bf{2967.8} & \textbf{2965.4} & \textbf{2966.5}\\

 \hline \hline
 
 cm201 & 5      & 4436.6 & 4452.5 & \bf{4367.6} & \bf{4400.1} & \textbf{4361.9} & \textbf{4374.4}\\
 cm202 & 6      & 4998.8 & 5024.9 & \bf{4983.2} & \bf{4990.0} & \textbf{4982.4} & \textbf{4987.4}\\
 cm203 & 5      & \bf{4445.8} & 4484.6& 4453.7  & \bf{4473.4} & \textbf{4440.3} & \textbf{4454.4}\\ 
 cm204 & 5      & 4335.2 & 4372.4 & \bf{4319.6} & \bf{4322.3} & \textbf{4317.5} & \textbf{4317.6}\\
 cm205 & 4      & 3863.5 & 3883.2 & \bf{3821.0} & \bf{3852.2} & \textbf{3816.8} & \textbf{3833.9}\\
 cm206 & 4      & 3722.0 & 3743.2 & \bf{3708.7} & \bf{3728.8} & \textbf{3699.0} & \textbf{3712.0}\\
 cm207 & 4      & 3968.4 & 3977.8 & \bf{3932.0} & \bf{3968.7} & \textbf{3917.4} & \textbf{3934.4}\\
 cm208 & 4      & 3771.1 & 3793.2 & \bf{3714.6} & \bf{3726.6} & \textbf{3714.6} & \textbf{3719.5}\\
 
 \hline \hline
 
 rcm101 & 10    & \bf{3062.0}&\bf{3062.0}&\bf{3062.0}& 3062.9 & \textbf{3062.0} & 3062.2\\
 rcm102 & 10    & 3132.2 & 3142.3 & \bf{3110.8} & \bf{3122.4} & \textbf{3110.8} & \textbf{3113.8}\\
 rcm103 & 10    & 3152.9 & 3163.8 & \bf{3121.3} & \bf{3127.0} & \textbf{3121.3} & \textbf{3121.9}\\ 
 rcm104 & 10    & 3119.6 & 3134.6 & \bf{3111.0} & \bf{3114.8} & \textbf{3111.0} & \textbf{3111.9}\\
 rcm105 & 10    & 3187.9 & 3210.7 & \bf{3166.2} & \bf{3173.0} & \textbf{3165.3} & \textbf{3169.5}\\
 rcm106 & 10    & 3218.9 & 3218.9 & \bf{3165.2} & \bf{3175.1} & \textbf{3158.2} & \textbf{3164.6}\\
 rcm107 & 11    & 3488.9 & 3514.0 & \bf{3487.7} & \bf{3494.4} & \textbf{3487.7} & \textbf{3490.2}\\
 rcm108 & 11    & 3592.7 & 3592.7 & \bf{3531.1} & \bf{3534.6} & \textbf{3531.1} & \textbf{3531.1}\\

 \hline \hline
 
 rcm201 & 2     & 2804.0 & 2827.8 & \bf{2738.8} & \bf{2758.8} & \textbf{2715.4} & \textbf{2748.1}\\
 rcm202 & 2     & 2836.9 & 2846.8 & \bf{2734.0} & \bf{2757.7} & \textbf{2716.8} & \textbf{2733.8}\\
 rcm203 & 2     & 2721.9 & \bf{2725.4} & \bf{2705.0} & 2739.9 & \textbf{2704.8} & \textbf{2714.2}\\ 
 rcm204 & 2     & 2726.5 & 2743.1 & \bf{2698.6} & \bf{2712.7} & \textbf{2692.6} & \textbf{2699.3}\\
 rcm205 & 2     & 2754.5 & 2775.7 & \bf{2729.5} & \bf{2748.3} & \textbf{2718.8} & \textbf{2725.1}\\
 rcm206 & 2     & 2812.7 & 2830.6 & \bf{2732.3} & \bf{2766.6} & \textbf{2721.9} & \textbf{2743.1}\\
 rcm207 & \bf{2}& 3749.8 & 3786.8 & \bf{2861.0} & \textbf{\bf{3049.9}} & \textbf{2861.0} & 3216.9\\
 rcm208 & 2     & 2791.4 & 2817.2 & \bf{2726.8} & \bf{2738.3} & \textbf{2722.7} & \textbf{2732.5}\\
 
 \hline \hline

\end{tabular}
\caption{Results on benchmark instances from \citep{Belhaiza2014AWindows} without time minimisation. The best known objectives (Best R.) using $\leq$ 60 seconds or best know solutions overall, are marked in bold. If our solver uses one fewer routes ($m$), the amount used is bold. All tests were repeated 10 times.}
\label{tab:resultsNTM}
\end{table}

 \begin{table}[htbp]
 \centering
 \footnotesize
 \begin{tabular}{c c | r r | r r | r r} 
   &  & \multicolumn{2}{c|}{\cite{Belhaiza2014AWindows}}  & \multicolumn{4}{c}{ALNS} \\
 \hline
 Instance & m & $B=1$  & & $B=1$  & 60 sec. & $B=1$  & 600 sec.\\
  &  & Best R. & Avg.  & Best R. & Avg. & Best R. & Avg. \\
 \hline
 
 rm101 & 10 & 4041.9 & 4072.9       & \bf{4014.2} & \bf{4035.3}             & \textbf{4014.0} & \textbf{4022.3} \\
 rm102 & 9  & 3765.1 & \bf{3765.1}  & \bf{3732.0} & 3771.1                  & \textbf{3729.9} & \textbf{3736.7} \\
 rm103 & 9  & 3708.5 & 3736.5       & \textbf{\bf{3699.7}} & \bf{3705.2}    & 3700.6 & \textbf{3701.5} \\ 
 rm104 & 9  & 3718.0 & 3722.9       & \textbf{\bf{3700.6}} & \bf{3702.0}    & \textbf{3700.6} & \textbf{3701.9} \\
 rm105 & 9  & 3688.8 & 3716.9       & \bf{3686.6} & \bf{3688.6}             & \textbf{3686.6} & \textbf{3687.0} \\
 rm106 & 9  & \textbf{\bf{3692.9}} & 3743.7  & 3708.0      & \bf{3713.3}    & 3708.0 & \textbf{3711.4} \\
 rm107 & 9  & 3701.4 & 3714.1       & \textbf{\bf{3689.9}} & \bf{3692.0}  & \textbf{3689.9} & \textbf{3689.9} \\
 rm108 & 9  & 3729.1 & 3729.1       & \textbf{\bf{3719.9}} & \bf{3725.6}  & \textbf{3719.9} & \textbf{3720.5} \\
 \hline \hline

 rm201 & \bf{2} & 4808.2 & 4847.4       & \bf{3885.6} & \bf{4089.0}  & \textbf{3815.0} & \textbf{3861.7} \\
 rm202 & 2      & \bf{3739.0} & 3775.7  & 3740.7      & \bf{3764.0}  & \textbf{3725.9} & \textbf{3742.2} \\
 rm203 & 2      & 3710.3 & 3728.3       & \bf{3696.9} & \bf{3721.5}  & \textbf{3696.1} & \textbf{3705.2} \\ 
 rm204 & 2      & \bf{3691.9} & 3708.8  & 3693.4 & \bf{3699.6}       & \textbf{3681.5} & \textbf{3687.8} \\
 rm205 & 2      & 3689.9 & 3707.7       & \bf{3679.8} & \bf{3691.7}  & \textbf{3678.7} & \textbf{3681.7} \\
 rm206 & 2      & 3703.4 & 3720.3       & \bf{3695.3} & \bf{3713.4}  & \textbf{3690.2} & \textbf{3697.6} \\
 rm207 & 2      & 3701.7 & 3719.4       & \bf{3691.8} & \bf{3705.5}  & \textbf{3688.3} & \textbf{3692.5} \\
 rm208 & 2      & 3682.8 & 3699.6       & \bf{3678.1} & \bf{3683.8}  & \textbf{3676.2} & \textbf{3682.6} \\
 \hline \hline

 cm101 & 10     & 12320.0     & \bf{12344.4}        & \bf{12270.4}& 12366.8       & \textbf{12253.2} & \textbf{12295.5} \\
 cm102 & 12     & \bf{12492.1}&\textbf{\bf{12492.1}}& 12554.8     & 12591.5       & \textbf{12467.5} & 12545.4 \\
 cm103 & \bf{11}& 12641.2     & 12687.7             & \bf{12568.5}& \bf{12655.1}  & \textbf{12506.0} & \textbf{12562.6} \\ 
 cm104 & \bf{13}& 13087.8     & 13117.9             & \bf{12913.0}& \bf{12969.1}  & \textbf{12911.9} & \textbf{12917.6} \\
 cm105 & \bf{10}& \bf{12083.4}&\textbf{\bf{12144.4}}& 12101.4     & 12144.5       & \textbf{12038.1} & \textbf{12064.6} \\
 cm106 & 10     & \bf{12073.9}& \bf{12133.9}        & 12096.9     & 12146.4  & \textbf{12043.6} & \textbf{12073.8} \\
 cm107 & \bf{10}& 12324.2     & 12364.1             & \textbf{\bf{12106.2}} & \bf{12147.9}  & \textbf{12106.2} & \textbf{12106.5} \\
 cm108 & 10     & 11990.4     & 12012.6             & \textbf{\bf{11984.3}}& \bf{11986.1}  & \textbf{11984.3} & \textbf{11984.3} \\
 \hline \hline
 
 cm201 & 5 & \bf{13520.1} & \bf{13591.7}& 13561.8 & 13654.0              & \textbf{13472.0} & \textbf{13551.2} \\
 cm202 & 6 & \bf{14027.3} & 14060.7     & 14045.7      & \bf{14060.5}    & \textbf{14021.3} & \textbf{14049.3} \\
 cm203 & 5 & 13497.2      & \bf{13512.8}& \bf{13489.8} & 13566.1         & \textbf{13465.5} & \textbf{13504.9} \\ 
 cm204 & 5 & 13359.8 & 13413.7          & \bf{13346.8} & \bf{13375.5}    & \textbf{13323.9} & \textbf{13342.3} \\
 cm205 & 4 & \bf{12884.1} & \bf{12963.1}& 12931.6      & 13018.7         & \textbf{12852.1} & \textbf{12914.9} \\
 cm206 & 4 & 12767.7 & 12811.2          & \bf{12753.4} & \bf{12809.1}    & \textbf{12729.8} & \textbf{12774.3} \\
 cm207 & 4 & 13009.7 & \bf{13017.6}     & \bf{12994.1} & 13060.7         & \textbf{12937.0} & \textbf{13007.2} \\
 cm208 & 4 & 12788.1 & 12805.2          & \bf{12737.7} & \bf{12764.8}    & \textbf{12729.8} & \textbf{12741.2} \\
 \hline \hline
 
 rcm101 & 10 & 4098.9 & 4129.7 & \textit{\bf{4079.6}} & \bf{4082.6}  & \textbf{4079.6} & \textbf{4081.0} \\
 rcm102 & 10 & 4222.6 & 4228.4 & \bf{4183.5} & \bf{4188.6}  & \textbf{4178.4} & \textbf{4180.8} \\
 rcm103 & 10 & 4174.3 & 4185.4 & \bf{4144.6} & \bf{4149.5}  & \textbf{4144.6} & \textbf{4146.5} \\ 
 rcm104 & 10 & 4156.3 & 4170.7 & \bf{4141.5} & \bf{4147.6}  & \textbf{4137.0} & \textbf{4141.0} \\
 rcm105 & 10 & 4216.7 & 4227.0 & \textbf{\bf{4188.4}} & \bf{4201.9}  & \textbf{4188.4} & \textbf{4195.0} \\
 rcm106 & 10 & 4219.9 & 4236.3 & \textbf{\bf{4175.4}} & \bf{4192.4}  & \textbf{4175.4} & \textbf{4179.4} \\
 rcm107 & 11 & 4542.4 & 4560.8 & \textbf{\bf{4516.5}} & \bf{4520.5}  & \textbf{4516.5} & \textbf{4517.0} \\
 rcm108 & 11 & 4614.5 & 4614.5 & \textbf{\bf{4565.2}} & \bf{4581.3}  & \textbf{4565.2} & \textbf{4569.2} \\
 \hline \hline
 
 rcm201 & 2     & \bf{3783.6} &\textbf{\bf{3824.5}}& 3816.0 & 4222.4      & \textbf{3733.8} & 3920.4 \\
 rcm202 & 2     & 3847.1 & \bf{3847.1}         & \bf{3792.6} & 4010.1      & \textbf{3756.1} & \textbf{3812.8} \\
 rcm203 & 2     &\bf{3721.9}&\textbf{\bf{3725.4}}& 3746.7      & 3810.4      & \textbf{3716.5} & 3779.7 \\ 
 rcm204 & 2     & 3726.5 & 3743.1              & \bf{3700.8} & \bf{3721.0} & \textbf{3699.8} & \textbf{3708.1} \\
 rcm205 & 2     & 3754.5 & 3775.7              & \bf{3732.0} & \textbf{\bf{3753.0}} & \textbf{3731.8} & 3756.6 \\
 rcm206 & 2     & 3812.7 & 3830.6              & \bf{3741.5} & \bf{3777.4} & 3744.7 & \textbf{3777.2} \\
 rcm207 & \bf{2}& 4764.2 & 4792.2              & \bf{3859.7} & \bf{4581.2} & \textbf{3859.1} & \textbf{4226.5} \\
 rcm208 & 2     & 3791.4 & 3817.2              & \textbf{\bf{3731.7}} & \textbf{\bf{3746.2}} & \textbf{3731.7} & 3747.5 \\
 \hline \hline

\end{tabular}
\caption{Results on benchmark instances from \citep{Belhaiza2014AWindows} using time minimisation, best known objectives (Best R.) using $\leq 60$ seconds or best overall are marked in bold. If our solver uses one fewer routes ($m$), the amount used is bold. All tests were repeated 10 times.}
 \label{tab:resultsTM}
\end{table}

The results in Table \ref{tab:resultsNTM} indicate that our \ac{ALNS} generally outperforms current state-of-the-art for VRPMTW without time minimisation, and provides new currently best known solutions for 40 of the 48 benchmark instances within 60 seconds. If the time limit is extended to 600 seconds, we find best known solutions for an additional 5 instances, totalling 45 of the 48 provided. Most of the solutions from the sets \textit{rm1xx} and \textit{rcm1xx} are identical to those obtained within 60 seconds, however, for the remaining sets the solutions are typically improved.

Table \ref{tab:resultsTM} shows that \ac{ALNS} generally outperforms current state-of-the-art also for the VRPMTW with time minimisation. Improved upper bounds are found in 36 of 48 cases for the average objective, and new currently best known solutions are found for 37 of the 48 instances. The gap between the average and the best found solutions, indicate a variance in the performance of the algorithm and a potential for further improvement. To further investigate this, the heuristic was also run with a 10 minute time limit using the parameters configuration from Table~\ref{tab:defaultConf}. As the solutions generally improve with the additional time, we can conclude that our hypothesis that further improvements are possible is correct. 
This resulted in new best known solutions for 47 of 48 instances.

In both problems, for six of the instances, ALNS also provides a solution with fewer routes than the previously best known, and five of these are new best known solutions.

\subsection{On the challenges inherent to the problem}

The choice of a time window $p\in P_i$ for visit $i\in V$ is a disjunction. Models that accommodate such constraints often use a big-M approach as in Constraint~\eqref{eq:waiting}. This modelling technique is typically expected to make mathematical models harder to solve, in part, due to a degeneration of the lower bound. This is also supported by the current model being incapable of providing good bounds for any relevant instance size.

Heuristics
effectiveness depends on their ability to navigate in the search space. Denote by $sol_{opt}$ the optimal solution and assume we are given an objective function $cost(sol)$ and a neighbourhood $neighbour(sol_i)$ where $sol_j\in neighbour(sol_i)$ if and only if $sol_j$ is reachable from solution $sol_i$ in the local search algorithm. If there exists a sequence (search path) $sol_0, sol_1 \dots sol_{n}$ such that $sol_{n}=sol_{opt}$ and $sol_i \in neighbourhood(sol_{i-1})$ $\forall i\in 1\dots n$ for any initial solution $sol_0$, we call the solution space connected under the neighbourhood operator.

As the number of vehicles is not limited, it is clear that the search space is connected under neighbourhood operators such as move, where a single visit is moved at a time. Each visit can be moved to its own vehicle where after the optimal solution can be assembled visit by visit. The disjunctions will thus not impose a disconnected search space in traditional local search algorithms.

Any local search progresses through a search path, but as the potential number of search paths is prohibitive, the search is heuristically guided through a partial exploration. If a search path $sol_0, sol_1 \dots sol_{n}$ has the property $cost(sol_0) > cost(sol_1)> \dots >cost(sol_{n})$ in a minimisation problem, the search will not encounter a local minimum along its search path. If no such search path exists, any heuristic looking for the optimum solution must accept worsening solutions during the search as part of its diversification strategy. The difficulty of finding good solutions is typically considered to be positively correlated with number of consecutive non-improving neighbours that must be visited during a search.

The VRPMTW exhibits a potential for an increased number of non-improving moves between good solutions. Given a partial solution where a subset of visits are placed consecutively due to geographical closeness, a traditional VRPTW instance would contain time windows causing visits to be placed into the only feasible part of the route, eg. early, middle or late in the route. In VRPMTW they may not be scheduled to be visited during the optimal time window, and changing this time window directly may not be feasible with respect to the rest of the route. Moving the visits one by one is likely to increase the cost of the solution significantly before they potentially are rescheduled during an alternative time window.

The problem of swapping groups of customers between time windows is further compounded by the fact that scheduling a visit and choosing its time window immediately, already imposes precedence and compatibility constraints with unscheduled and thus unknown visits, as their time windows may not be compatible with promising orderings of the visits anymore. Our attempt to remedy these issues proved  important for both the average cost of solutions (Figure\ref{fig:omit}) and in the search for new best solutions (Figire \ref{fig:omitBest}).

\section{Conclusions}
\label{sec:conclusions}
We have formulated a more compact mathematical model for both variants of the \ac{VRPMTW}, though bounds and solutions remain unobtainable the model is more readable, and an error regarding service times has been corrected.

We have presented an adaptation of \ac{ALNS} to two variants of the VRPMTW. We have shown that the problem without time minimisation effectively can be solved using traditional heuristics based on insertion with an overhead scaling linearly in the number of time windows. We have furthermore found new best solutions for all but 6 instances.

We presented a labelling algorithm for evaluating delta costs of insertions for VRPMTW with time minimisation. Furthermore we demonstrated that \ac{ALNS} using the delta evaluation is competitive with existing methods from the literature, and for all categories of instances outperform them. Lastly we found new best known solutions to all but 2 instances.

We have observed that an increased focus on diversification over intensification is needed when tuning the algorithm. Furthermore we noticed a large variance in the quality of results, making us suspect further work is warranted. This suspicion was confirmed when running the algorithm for 10 minutes, confirming that better solutions are obtainable.

\section*{Acknowledgements}

This project was partially funded by the Danish Transport Authority under the grant TS20707-00132.

\section*{References}
\bibliography{references}

\end{document}